\documentclass[psamsfonts,11pt]{amsart}

\usepackage{amssymb,amsfonts, amsmath}
\usepackage[all,arc]{xy}
\usepackage{enumerate}
\usepackage{mathrsfs}
 
\usepackage[margin=1.09in,footskip=0.5in]{geometry}
\newtheorem{thm}{Theorem}[section]
\newtheorem{cor}[thm]{Corollary}
\newtheorem{prop}[thm]{Proposition}
\newtheorem{lem}[thm]{Lemma}

\newtheorem*{thm:associativity}{Theorem \ref{thm:associativity}}

\theoremstyle{definition}

\newtheorem{defns}[thm]{Definitions}

\newtheorem{exmp}[thm]{Example}

\theoremstyle{remark}
\newtheorem{rem}[thm]{Remark}
\newtheorem{rems}[thm]{Remarks}

\makeatletter
\let\c@equation\c@thm
\makeatother
\numberwithin{equation}{section}

\renewcommand{\setminus}{\smallsetminus}

\title[Euler Characteristics on virtually free products]{Euler Characteristics on virtually free products}

\subjclass[2010]{20E06, 20E26.}

\keywords{Euler characteristics, residually finite groups, virtually torsion free groups, rank gradient, first ${L^2-}$betti number.}

\author[Konstantinos Tsouvalas]{\bfseries Konstantinos Tsouvalas}

\address{
Department of Mathematics, National and Kapodistrian University of Athens, Panepistimiopolis, Athens 15784, Greece}
\email{\texttt{kostastsouvalas@windowslive.com}}

\date{\today}

\begin{document}

\begin{abstract} We define Euler characteristics on classes of residually finite and virtually torsion free groups and we show that they satisfy certain formulas in the case of amalgamated free products and HNN extensions over finite subgroups. These formulas are obtained from a general result which applies to the rank gradient and the first ${L^2-}$Betti number of a finitely generated group.

\end{abstract}

\maketitle

\section{Introduction}
Let ${\mathcal{C}}$ be a class of groups closed under subgroups of finite index. An Euler characteristic on ${\mathcal{C}}$ is a mapping ${\chi:\mathcal{C}\rightarrow R}$, where ${R}$ is a commutative ring with identity, such that if ${G}$ is a group in ${\mathcal{C}}$ and ${H}$ is a finite index subgroup of ${G}$, then ${\chi\left ( H \right )=\left [ G:H \right ]\chi\left ( G \right )}$ and ${\chi\left (G_{1} \right )=\chi\left ( G_{2} \right )}$ if ${G_{1}, G_{2}}$ are isomorphic groups in ${\mathcal{C}}$. An Euler characteristic with values on a field, offers information about the isomorphism classes of subgroups: if ${\chi\left ( G \right )\neq 0}$ and ${H_{1}, H_{2}}$ are two isomorphic and finite index subgroups of ${G}$, then they have the same index. In this case, every monomorphism ${\varphi:G\rightarrow G}$ is either an automorphism or the image of ${\varphi}$ has infinite index in ${G}$. Euler characteristics have been defined on various classes of groups, satisfying some homological conditions, by Brown ${[1]}$, Serre ${[13]}$, Stallings ${[15]}$ and Wall ${[18]}$ (see also ${[2]}$ and the references therein). In some cases, they satisfy the following properties for amalgamated free products and HNN extensions: If ${A,B \in \mathcal{C}}$ and ${H}$ is a finite group, then

\vspace{0.2cm} 

${\textup{(1)}}$ ${A \ast_{H} \in \mathcal{C}}$ and ${\chi\left ( A\ast_{H} \right )={\chi}\left (A \right )- \frac{1}{|H|}}$.

\vspace{0.2cm}

${\textup{(2)}}$ ${A \ast_{H} B \in \mathcal{C}}$ and ${\chi\left ( A\ast_{H} B \right )={\chi}\left (A \right )+{\chi}\left (B \right )-\frac{1}{|H|}}$. 

\vspace{0.2cm}

In this paper, following the point of view of ${[11]}$ and ${[16]}$, we define Euler characteristics on classes of residually finite or virtually torsion free groups which are non-zero on each finitely generated residually finite group ${G}$ with infinitely many ends. In particular, we conclude that isomorphic subgroups of finite index of ${G}$ have the same index.
Moreover, we prove that they also satisfy the above properties. This is deduced as a special case of the following theorem:

\begin{thm:associativity}

Let ${\mathcal{C}}$ be a class of groups closed under free products and subgroups of finite index. We assume that ${\mathcal{C}}$ contains the infinite cyclic group ${\mathbb{Z}}$. Suppose that ${\chi:\mathcal{C}\rightarrow \mathbb{R}}$ is an Euler characteristic which satisfies the following property: ${\chi\left ( A \ast B\right )=\chi\left ( A\right )+\chi\left ( B\right )-1}$ for every ${A,B \in \mathcal{C}}$.

\vspace{0.2cm}

${\textup{(i)}}$ Let ${G\ast_{H}}$ be an HNN extension, where ${G}$ is residually finite or a virtually torsion free group in ${\mathcal{C}}$ and ${H}$ finite. Then ${G\ast_{H}}$ is virtually in ${\mathcal{C}}$ and: $${\chi\left ( G\ast_{H} \right )={\chi}\left ( G \right )- \frac{1}{|H|}}$$

\vspace{0.1cm}
${\textup{(ii)}}$ Let ${G_{1}\ast_{H} G_{2}}$ be an amalgamated free product, where ${G_{1}, G_{2}}$ are residually finite or virtually torsion free groups in ${{\mathcal{C}}}$ and ${H}$ finite. Then ${G_{1}\ast_{H} G_{2}}$ is virtually in ${\mathcal{C}}$ and: $${\chi\left ( G_{1} \ast_{H} G_{2} \right )={\chi}\left ( G_{1} \right )+{\chi}\left ( G_{2} \right )- \frac{1}{|H|}}$$

\end{thm:associativity}

\vspace{0.2cm}

\noindent Note that in the above theorem, we extend the Euler characteristic $\chi$ on virtually ${\mathcal{C}}$ groups, by using Wall's argument. More precisely, if $G$ is a group containing a finite-index subgroup ${K}$ in ${\mathcal{C}}$, then we define ${\chi\left ( G \right )=\frac{\chi\left (K \right )}{\left [ G:K \right ]}}$. \par

We show that the above theorem applies to the rank gradient of a group ${[4]}$  and we use it to calculate the rank gradient of the fundamental group of a finite graph of groups with residually finite or virtually torsion free vertex groups and finite edge groups. Moreover, the main theorem applies to the first ${L^2-}$betti number ${[6,9]}$ of a finitely generated group. By using the formula for amalgamated free products in ${[7]}$, we deduce the corresponding formulas for HNN extensions over finite subgroups.

As another application, we obtain a proof of Linnell's accessibility result ${[5]}$ in the case of residually finite and virtually torsion free groups. We also prove in some cases the accessibility of automorphism fixed subgroups.

\section{Finite-index subgroups of amalgamated free products and HNN extensions}
First, we recall some basic facts and definitions on classifying spaces of groups. For a group ${G}$, we denote by ${EG}$ the topological space defined as follows: ${EG}$ is the contractible simplicial complex whose ${n-}$simplices are ordered ${n-}$tuples ${\left  [g_{0},...,g_{n} \right ]}$ of elements of ${G}$, with the natural inclusions. The group ${G}$ acts simplicially on the space ${EG}$ by left multiplication: ${g\cdot \left [ g_{0},...,g_{n} \right ]=\left [gg_{0},...,gg_{n}\right ]}$ and we denote by ${BG}$ the quotient space ${EG/G}$. The action is free and each simplex has trivial stabilizer. Any subgroup ${H}$ of ${G}$ acts on ${EG}$. Let ${p_{H}:EG\rightarrow  EG/H}$ denote the corresponding quotient map. Note also that ${EH}$ and ${BH}$ can be viewed as subspaces of ${EG}$, ${BG}$ respectively. For any groups ${H_{1}, H_{2}}$ and each homomorphism ${\phi: H_{1}\rightarrow H_{2}}$, there is the induced map ${E\phi: EH_{1} \rightarrow EH_{2}}$ which sends the simplex ${\left [ h_{0},...,h_{n} \right ]}$ linearly onto the simplex ${\left [ \phi(h_{0}),...,\phi(c_{n}) \right ]}$. In the case where ${K}$ is a normal subgroup of ${G}$, the quotient group ${G/K}$ acts on the quotient space ${EG/K}$ in the obvious way. We denote by ${p_{G/K}}$ the corresponding quotient map. We refer the reader to ${[3]}$ for further details.

\vspace{0.2cm}

The following proposition is essential for the computations of the fundamental groups in the following lemmas:

\vspace{0.2cm}

\begin{prop}\emph{}Let ${K,H}$ be two subgroups of a group ${G}$ such that ${K}$ is normal and ${H \cap K=1}$. Then the inverse image ${p_{G/K}^{-1}\left ( p_{G}\left ( EH \right ) \right )}$ is a disjoint union of contractible subspaces of ${p_{K}\left ( EG \right )}$.\end{prop}

\vspace{0.1cm}

\begin{proof} We first note that ${H}$ acts on ${G/K}$ with left multiplication. Thus there is a partition of ${G/K}$ into ${H-}$orbits (double cosets); i.e. ${G/K =\bigsqcup_{i \in I} H g_{i}K}$. Since ${H \cap K=1}$, the subspace ${EH}$ contains at most one point of each ${K-}$orbit and therefore  the restriction ${p_{K}|_{EH}}$ is a homomorphism onto its image, since ${p_{K}}$ is a covering map. In particular, ${p_{K}\left ( EH \right )}$ is contractible as well as each translations ${g_{j}K  \cdot p_{K}\left ( EH \right )}$. We now suppose that the intersection of the translations ${g_{i}K \cdot p_{K}\left ( EH \right ) \cap g_{j}K \cdot p_{K}\left ( EH \right )}$ is non-empty. This means that there are elements ${h_{1},...,h_{m},h_{1}',...,h_{d}' }$ in ${H}$ and an element ${k}$ of ${K}$ such that the intersection of the simplices ${\left [ kg_{i}h_{1},...,kg_{i}h_{m}\right ]}$, ${\left [g_{j}h_{1}',...,g_{j}h_{d}'\right ]}$ is non-empty. It follows that the elements ${g_{i}}$ and ${g_{j}}$ define the same double coset and hence $${p_{G/K}^{-1}\left ( p_{G}\left ( EH \right ) \right )=\bigsqcup_{i \in I}  g_{i}K  \cdot p_{K}\left ( EH \right ).}$$ \end{proof}

\vspace{0.2cm}

Now we will construct suitable covering maps in order to obtain finite index subgroups of amalgamated free products and HNN extensions over finite groups.

\vspace{0.2cm}

\begin{lem}\emph{•}Let ${G_{1}, G_{2}}$ be groups with subgroups of finite index ${H_{1}, H_{2}}$, respectively with ${m=\left [ G_{1}:H_{1} \right ],n=\left [ G_{2}:H_{2} \right ]}$. Then, for any common multiple ${s}$ of ${m,n}$, there exists a subgroup ${K}$ of the free product ${G_{1}\ast G_{2}}$ such that: $${K\cong \underbrace{H_{1}\ast...\ast H_{1}}_{s/m}\ast\underbrace{H_{2}\ast...\ast H_{2}}_{s/n}\ast F_{d}}$$ where ${d=\left ( \frac{s}{m}-1 \right )\left ( \frac{s}{n}-1 \right )}$ and ${\left [ G_{1}\ast G_{2}:K \right ]=s}$. \end{lem}

\vspace{0.1cm}

\begin{proof} Let ${X_{1}, X_{2}}$ be two CW complexes such that ${G_{1}=\pi_{1}\left ( X_{1},x_{1} \right )}$ and  ${G_{2}}=$ ${\pi_{1}\left ( X_{2},x_{2} \right )}$. We consider the space which is obtained from the disjoint union ${X_{1}\sqcup X_{2} \sqcup \left [ 0,1 \right ]}$ via the identifications ${0 \sim x_{1}, 1 \sim x_{2}}$. For each subgroup ${H_{i}}$, let ${p_{i}:X_{H_{i}}\rightarrow X_{i}}$ be the corresponding covering map such that ${p_{i \ast}\left ( \pi_{1}\left ( X_{H} \right ) \right ) \cong H_{i}}$. Let ${Z_{1}=\sqcup_{i=1}^{s/m}X_{H_{1},i}}$ and ${Z_{2}= \sqcup_{j=1}^{s/n}X_{H_{2},j}}$ be the disjoint unions of ${s/m}$ and ${s/n}$ copies of the spaces ${X_{H_{1}}}$ and ${X_{H_{2}}}$ respectively. Let also ${x_{1r}^{(i)}}$ and ${x_{2t}^{(j)}}$ be elements of the copies ${X_{H_{1},i}}$ and ${X_{H_{2},j}}$ respectively such that ${p_{1}\left ( x_{1r}^{\left ( i \right )} \right )=x_{1},p_{2}\left ( x_{2t}^{\left ( j \right )} \right )=x_{2}}$ where $1\leqslant r\leqslant m,$ $1\leqslant t \leqslant n$. By using arcs ${K_{1},...K_{s}}$, we construct a path connected covering space ${\widetilde{X}}$ of ${X}$ with exactly ${s}$ sheets with the following way: the arc ${K_{ij}}$ identifies the ${j-}$th element ${x_{1j}^{(i)}}$ of the ${i-}$copy of ${X_{H_{1}}}$ with the ${i-}$th element ${x_{2i}^{(j)}}$ of the ${j-}$copy of ${X_{H_{2}}}$. Let ${f_{ij}:K_{ij}\rightarrow \left [ x_{1},x_{2} \right ]}$ be homeomorphisms with ${f_{ij}\left (x_{1i}^{(j)} \right )=x_{1},f_{ij}\left (x_{2j}^{(i)} \right )=x_{2}}$. Then ${\widetilde{X}}$ is a ${s-}$sheeted covering space via the covering map ${p:\widetilde{X} \rightarrow X}$ defined as follows: ${p|_{X_{H_{1},i}}=p_{1},\ p|_{X_{H_{2},j}}}$ ${=p_{2},\ p|_{K_{i}}=f_{i}}$. Finally, by using the Seifert-van Kampen theorem we have: $${\pi_{1} ( \widetilde{X})\cong \ast_{i=1}^{s/m}H_{1} \ast_{i=1}^{s/n}H_{2}\ast F_{d}}$$ where ${d=\left ( \frac{s}{m}-1 \right )\left ( \frac{s}{n}-1 \right )}$ and the conclusion follows.\end{proof}

\vspace{0.2cm}

\begin{lem}\emph{}Let ${G}$ be a group, ${C}$ a finite subgroup of ${G}$, ${\phi: C\rightarrow G}$ a monomorphism and ${K}$ a normal and finite index subgroup of ${G}$ such that ${K \cap C=K \cap\phi(C)=1}$. Then the HNN extension ${G\ast_{C}}$ relative to ${\phi}$ has a subgroup ${N}$ of index ${\left [ G:K \right ]}$ which is isomorphic to the free product ${K \ast F_{\frac{\left [ G:K \right ]}{|C|}}}$. \end{lem}

\vspace{0.1cm}

\begin{proof} We consider the space $${X=BG \sqcup \big(T \times [0,1]\big) \big/ \big\{\left ( \left [x  \right ]_{G},0 \right )\sim \left [x  \right ]_{G}, \ \left ( \left [x  \right ]_{G},1 \right )\sim \left [ E\phi\left ( x \right ) \right ]_{G} \big\} }$$ where ${T=p_{G}\left ( EC \right )}$, ${S=p_{G}\left ( E\phi(C) \right )}$. The space ${p_{K}\left ( EG \right )}$ is a ${\left [ G:K \right ]-}$sheeted covering space of ${BG}$ with covering map ${p_{G/K}}$. By Proposition 2.1, we have the partitions $${p_{G/K}^{-1}\left ( T \right )=\bigsqcup_{i \in I}g_{i}K \cdot p_{K}\left ( EC \right ),\ p_{G/K}^{-1}\left ( S \right )=\bigsqcup_{i \in I}h_{i}K \cdot p_{K}\left ( E\phi(C) \right )}$$ where ${|I|=\frac{[G:K]}{|C|}}$, ${K \cap C=K \cap\phi(C)=1}$. We glue the cylinders $${L_{i} \times\left [ 0,1 \right ], \ L_{i}=g_{i}K \cdot p_{K}\left ( EC \right )}$$ to the space ${p_{K}\left ( EG \right )}$ by making the following identifications: $${\Big( \left (g_{i}K  \right )\cdot \left [ x \right ]_{K},0 \Big)\sim \left (g_{i}K  \right )\cdot \left [ x \right ]_{K}, \ \Big( \left (g_{i}K  \right )\cdot \left [ x \right ]_{K},1 \Big )\sim \left (h_{i}K  \right )\cdot \left [ E \phi\left (x  \right ) \right ]_{K}}$$ where ${p_{K}\left ( x \right )=\left [ x \right ]_{K}}$. Let ${\widetilde{X}}$ be the resulting quotient space. \par
\noindent Now let ${q_{1}: BG \sqcup \left (T \times [0,1]  \right )\rightarrow X}$, ${q_{2}:p_{K}\left ( EG \right ) \sqcup \bigsqcup_{i \in I} L_{i}\times [0,1] \rightarrow \widetilde{X}}$ be the  natural quotient maps and let ${p:\widetilde{X} \rightarrow X}$ be the continuous map defined as follows:

\vspace{0.1cm}
\begin{center}
 ${p\left (q_{2}\left ( y \right ) \right )=\left\{\begin{matrix}q_{1}\left ( p_{G/K}\left ( s \right )  \right ), \  \textup{if} \  y=s \in p_{K}\left ( EG \right )  \ \ \ \ \ \ \ \\ 
 \\ q_{1}\left ( p_{G/K}\left ( s \right ),t  \right ), \  \textup{if} \  y=(s,t) \in L_{i}\times [0,1] \
\end{matrix}\right.}$
\end{center}

\vspace{0.1cm}
\noindent In order to prove that ${p}$ is a covering map, it is enough to consider the cases where an arbitrary element ${y}$ lies in one of the subspaces ${T, S, q_{1}\left ( T\times\left ( 0,1 \right ) \right )}$. Let ${y}$ be an element in the union ${S \cup T}$. Since ${p_{G/K}}$ is a covering map, there exists an elementary neighbourhood ${U_{y}}$ of ${y}$, i.e. ${p_{G/K}^{-1}\left (U_{y} \right )=\sqcup_{i \in J} V_{i}}$  where each ${V_{i}}$ is mapped homeomorphically onto ${U_{y}}$. In the case where ${y}$ be an element of ${T}$, we consider the set ${W_{y}=U_{y} \sqcup \left ( q_{1}\left ( U_{y}\cap T \times (0, \varepsilon) \right ) \right )}$, where ${0<\varepsilon<1}$. The set ${W_{y}}$ is open since its inverse image ${q_{1}^{-1}(W_{y})=U_{y} \sqcup \big(\left ( U_{y}\cap T \right )\times \left [ 0, \varepsilon \right )\big)}$ is open and ${q_{1}}$ is a quotient map. We also observe that the inverse image ${p^{-1}(W_{y})}$ is the disjoint union of the sets: $${W_{i}=V_{i} \sqcup q_{2}\left (\left (\cup_{j} V_{i}\cap {L_{j}} \right ) \times (0, \varepsilon) \right ),\ i \in J}$$ It is easy to see that the sets ${W_{i}}$ are homeomorphic with ${W_{y}}$ via ${p}$ and that ${q_{2}^{-1}(W_{i})=}$ ${V_{i} \sqcup  \left ( \left (\cup_{j} V_{i}\cap {L_{j}}\right ) \times [0, \varepsilon)\right )}$. It follows that the sets ${W_{i}}$ are open. In the case where ${y}$ lies in ${S}$, we consider the open region $${U_{y}\sqcup q_{1}\left ( \left (p_{G} \circ E\phi  \right )^{-1}\left ( U_{y}\cap S \right )\times (1-\varepsilon,1)\right )}$$ and work similarly. If ${y=q_{1}(x,t)}$ for some ${0<t<1}$, we consider an elementary neighbourhood ${U_{x}}$ of ${x}$ with ${p_{G/K}^{-1}\left ( U_{x} \right )=\sqcup_{i \in J} U_{i}}$. For ${\varepsilon>0}$ small enough, the set ${W_{x}=q_{1}\left ( U_{x} \cap T \times \left ( t-\varepsilon, t+\varepsilon \right )\right )}$ is an open neighbourhood of ${y}$ and the inverse image ${p^{-1}\left ( W_{x} \right )}$ is the disjoint union of the sets: $${E_{i}=q_{2}\left ( \left (\cup_{j} U_{i}\cap {L_{j}}\right )\times (t-\varepsilon, t+\varepsilon) \right ), \ i \in J}$$ Clearly the sets ${E_{i}}$ are open and homeomorphic (via the restriction) to ${W_{x}}$ and finally ${\widetilde{X}}$ is a ${\left [ G:K \right ]-}$sheeted covering space of ${X}$. Since the subspaces ${g_{i}K \cdot p_{K}\left ( EC \right ),}$ ${h_{i}K \cdot p_{K}\left ( E\phi\left ( C \right ) \right )}$ are contractible, the Seifert-van Kampen Theorem implies that ${\pi_{1} ( \widetilde{X} )}$ is isomorphic to the free product ${ K \ast F_{\frac{\left [ G:K \right ]}{|C|}}}$ and the conclusion follows. \end{proof}

\vspace{0.2cm}

\begin{lem}\emph{}Let ${G_{1}, G_{2}}$ be groups and monomorphisms ${\phi:C\rightarrow G_{1}, \psi:C\rightarrow G_{2}}$ where ${C}$ is a finite group. Let ${N_{1},N_{2}}$ be normal and finite index subgroups of ${G_{1},G_{2}}$ respectively such that ${N_{1}\cap \phi\left ( C \right )=1, N_{2}\cap \psi\left ( C \right )=1}$. Then the amalgamated free product ${G_{1} \ast_{C}G_{2}}$ with respect to ${\phi}$ and ${\psi}$ contains a subgroup of index ${\frac{n_{1}n_{2}}{c}}$ which is isomorphic to the group: $${\underbrace{N_{1}\ast ...\ast N_{1}}_{n_{2}/c} \ast \underbrace{N_{2}\ast ...\ast N_{2}}_{n_{1}/c}\ast F_{s}}$$ where ${n_{1}=\left [ G_{1}:N_{1} \right ],n_{2}=\left [ G_{2}:N_{2} \right ], c=|C|}$ and ${s=\left ( \frac{n_{1}}{c}-1 \right )\left ( \frac{n_{2}}{c}-1 \right )}$.\end{lem}

\vspace{0.1cm}

\begin{proof} Let ${X=BG_{1} \sqcup  BG_{2} \sqcup \big(BC \times [0,1] \big) \big/ \big\{\left ( \left [x  \right ]_{C},0 \right )\sim \left [E\phi\left (x   \right ) \right ]_{G_{1}}, \left ( \left [x  \right ]_{C},1 \right )\sim \left [ E\psi\left ( x \right ) \right ]_{G_{2}} \big\}  }$. We consider the disjoint unions $${X_{1}=\sqcup _{j=1}^{n_{2}/c}{\left ( p_{N_{1}}\left ( EG_{1} \right ) \right )_{j}}, \ X_{2}=\sqcup _{i=1}^{n_{1}/c} {\left (p_{N_{2}}\left ( EG_{2} \right ) \right )_{i}}}$$ of ${n_{1}/c}$ and ${n_{2}/c}$ copies of the spaces ${ p_{N_{1}}\left ( EG_{1} \right ),  p_{N_{2}}\left ( EG_{2} \right )}$ respectively. By Proposition 2.1 we have that $${p_{G_{1}/N_{1}}^{-1}\left ( p_{G_{1}}\left ( E\phi\left ( C \right ) \right ) \right )=\bigsqcup_{i \in I_{1}}g_{i}N_{1} \cdot p_{N_{1}}\left ( E\phi\left ( C \right ) \right )}$$ $${p_{G_{2}/N_{2}}^{-1}\left ( p_{G_{2}}\left ( E\psi\left ( C \right ) \right ) \right )=\bigsqcup_{j \in I_{2}}h_{j}N_{2} \cdot p_{N_{2}}\left ( E\psi\left ( C \right ) \right )}$$ where ${\left | I_{1} \right |=\frac{n_{1}}{c}, \left | I_{2} \right |=\frac{n_{2}}{c} }$. We also consider the cylinders $${C_{ij}=\left ( EC\times \left [ 0,1 \right ] \right )_{ij}}$$ where ${i \in I_{1}}$ and the index ${j}$ ranges over the copies of ${p_{N_{1}}\left ( EG \right )}$. Then the cylinders ${C_{ij}}$ are glued on the disjoint union ${X_{1}\sqcup X_{2}}$ with the following way: 
for every pair ${(i,j)}$ and element ${x \in EC}$, ${\left ( \left [ x \right ]_{C},0 \right )_{ij}}$ is identified with the element ${\left (g_{i}N_{1}\cdot p_{N_{1}}\left ( E \phi\left ( x \right ) \right )  \right )_{j}}$ as well as ${\left ( \left [ x \right ]_{C},1 \right )_{ij}}$ is identified with $\left (h_{j}N_{2}\cdot p_{N_{2}}\left ( E \psi\left ( x \right ) \right )  \right )_{i}$. Let ${\widetilde{X}}$ be the quotient space. It is clear from the previous construction that ${\widetilde{X}}$ is path connected and let ${p}$ be the map defined as follows: $${p\left (q_{2}(y) \right )=\left\{\begin{matrix}
q_{1}\left ( p_{G_{k}/N_{k}}(s) \right ),\ \textup{if} \ y=s \in X_{1}\sqcup X_{2}, \ k \in \left \{ 1,2 \right \} 
 \\
\\
q_{1}\left (\left [ s \right ]_{C},t \right ),  \  \textup{if} \ y=\left ( s,t \right )_{ij} \in C_{ij}, i \in I_{1}, j \in I_{2} \  \ \ \ \\ \ \ \ \ \
\end{matrix}\right.}$$ where ${q_{1}:BG_{1} \sqcup BG_{2} \sqcup \big(BC \times [0,1] \big) \rightarrow X}$, ${q_{2}:X_{1}\sqcup X_{2} \sqcup \sqcup_{ij} C_{ij}\rightarrow \widetilde{X}}$ are the natural quotient maps. The map ${p}$ is continuous by the gluing lemma and it is proved similarly as before that it is a covering map. Hence ${\widetilde{X}}$ is a ${\frac{n_{1}n_{2}}{c}-}$sheeted covering space of ${X}$. By Proposition 2.1, the translations ${g_{j}N_{1} \cdot p_{N_{1}}\left ( E \phi\left ( C \right ) \right ),}$ ${h_{i}N_{2} \cdot p_{N_{2}}\left ( E \psi\left ( C \right ) \right )}$ are contractible, hence the space ${\widetilde{X}}$ looks like a graph of spaces where the base graph has ${\frac{n_{1}n_{2}}{c^2}}$ edges and ${\frac{n_{1}+n_{2}}{c}}$ vertices. Thus, the Seifert-van Kampen theorem implies that the fundamental group ${\pi_{1} ( \widetilde{X})}$ is isomorphic to the free product $${\underbrace{N_{1}\ast ...\ast N_{1}}_{n_{2}/c} \ast \underbrace{N_{2}\ast ...\ast N_{2}}_{n_{1}/c}\ast F_{s}}$$ where ${s=\left ( \frac{n_{1}}{c}-1 \right )\left ( \frac{n_{2}}{c}-1 \right ).}$ \end{proof}

\section{Main Results}
Let ${G=\ast_{i \in I}G_{i}}$ be the free product of a family of groups ${\left \{ G_{i}: i \in I \right \}}$ and let ${H}$ be a subgroup of ${G}$. By Kurosh subgroup theorem, ${H}$ is a free product ${\ast_{j \in J}\left ( g_{j}H_{j}g_{j}^{-1} \right )\ast F\left ( X \right )}$, where ${H_{j}}$ is a subgroup of some free factor $G_{i}$ and ${ F\left ( X \right )}$ is a free group on ${X}$. In the case where the sets ${X, J}$ are finite, the Kurosh-rank of ${H}$ is defined as ${\textup{Kr}\left (H \right )=|X|+\left |J  \right |}$ where each ${H_{j}}$ is assumed to be non-trivial. If the number of factors $|I|$ is finite and ${H}$ is a finite-index subgroup of ${G}$, then, by ${[16, \textup{Prop. 3.2}]}$, the Kurosh-rank of ${H}$ with respect to the given splitting is ${1+\left [ G:H \right ]\left ( n-1 \right )}$.

\vspace{0.2cm}

For actions of groups on simplicial trees and the definition of the fundamental group of a graph of groups, we refer the reader to ${[14]}$. For a graph of groups ${\left ( \mathcal{G},Y \right )}$, we denote by ${G_{u}}$ the vertex group assigned to the vertex ${u \in VY}$ and by ${G_{e}}$ the assigned group to the edge ${e \in EY}$.

\vspace{0.2cm}

For a finitely generated group ${G}$, we denote by ${r(G)}$ the size of a minimal generating set of ${G}$.

\vspace{0.2cm}

\begin{defns} Let ${\mathcal{F}}$ be the class consisting of all groups which are free products with finitely many factors such that subgroups of finite index in each factor are freely indecomposable. Clearly, by Kurosh subgroup theorem ${\mathcal{F}}$ is closed under finite index subgroups. Examples of groups in ${\mathcal{F}}$ include finitely generated free groups and finitely generated torsion free groups. By ${\textup{V} \mathcal{F}}$ we denote the class of groups which are virtually  ${\mathcal{F}}$. \par Let ${H}$ be a group in ${\mathcal{F}}$ and ${H =\ast_{i=1}^{n} H_{i}}$ a splitting of ${H}$ as above. If ${H}$ is a subgroup of finite index of a group ${G}$, then we define $${\omega \left(G \right )=\frac{1-n}{\left [ G:H \right ]}}$$ Thus we have a mapping ${\omega: \textup{V} \mathcal{F}\rightarrow \mathbb{Q}}$. \par  

\end{defns}

\vspace{0.2cm}

\begin{defns} For a group ${G}$ we denote by ${\Lambda_{G}}$ the set of finite index subgroups directed by reverse inclusion. Let ${\mathcal{C}}$ be a class of groups and ${\sigma:\mathcal{C}\rightarrow \mathbb{R}}$ a mapping. We define the ${\sigma-}$volume to be the mapping ${V_{\sigma}:\mathcal{C}\rightarrow \overline{\mathbb{R}}}$ where $${V_{\sigma}\left ( G \right )=\varlimsup_{H \in \Lambda_{G}}\frac{\sigma\left ( H \right )}{\left [ G:H \right ]}=\inf_{H \in \Lambda_{G}} \sup_{K \in \Lambda_{H}}\frac{\sigma\left ( K \right )}{\left [ G:K \right ]}}$$ We also define the lower and the upper volumes respectively: $${\underline{V_{\sigma}}\left ( G \right )=\sup_{H \in \Lambda_{G}}\frac{\sigma\left ( H \right )}{\left [ G:H \right ]}, \ \overline{V_{\sigma}}\left ( G \right )=\inf_{H \in \Lambda_{G}}\frac{\sigma\left ( H \right )}{\left [ G:H \right ]}}$$ and we easily deduce that for every finite index subgroup ${H}$ of ${G}$ we have $${\underline{V_{\sigma}}\left ( H \right )\leqslant \left [ G:H \right ]\underline{V_{\sigma}}\left ( G\right ), \ \overline{V_{\sigma}}\left ( H \right )\geqslant  \left [ G:H \right ]\overline{V_{\sigma}}\left ( G \right )}$$ For a mapping ${\sigma: \mathcal{C}\rightarrow \mathbb{R}}$, we define the mapping ${\widetilde{V_{\sigma}}}$ on the class ${\mathcal{C}\cup\left \{ \textup{finite groups}\right \}}$: $${\widetilde{V_{\sigma}}\left ( G \right )=\left\{\begin{matrix}
V_{\sigma}\left ( G \right ),\ \textup{if} \ G \ \textup{is infinite in} \ \mathcal{C} \\ 
 -\frac{1}{\left | G \right |},\ \textup{if} \ G \ \textup{is finite} \ \ \ \ \ \ 
\end{matrix}\right.}$$
\end{defns}

\vspace{0.2cm}

\begin{prop}\emph{•} Let ${\mathcal{C}}$ be a class of groups closed under subgroups of finite index and ${\sigma: \mathcal{C}\rightarrow \mathbb{R}}$ a mapping such that for every ${G}$ in ${\mathcal{C}}$, the net ${\left (\frac{\sigma\left ( H \right )}{\left [ G:H \right ]}\right )_{H}}$ is bounded. Then ${V_{\sigma}}$ is an Euler characteristic.  \end{prop}
\vspace{0.1cm}
\begin{proof} By definition, we have ${V_{\sigma}\left ( G \right )=\inf_{L<G}\frac{1}{\left [ G:L \right ]}\underline{V_{\sigma}}\left ( L \right )}$. Let ${H}$ be a finite index subgroup of ${G}$ and ${\varepsilon>0}$. There exists a finite index subgroup ${L}$ of ${H}$ such that $${V_{\sigma}\left ( H \right )+\varepsilon\geqslant \frac{1}{\left [ H:L \right ]}\underline{V_{\sigma}}\left ( L \right )\geqslant \left [ G:H \right ]V_{\sigma}\left ( G \right )}$$ On the other hand, there exists a finite index subgroup ${K}$ of ${G}$ such that $${V_{\sigma}\left ( G \right )+\varepsilon\geqslant \frac{1}{\left [ G:K \right ]}\underline{V_{\sigma}}\left ( K \right )\geqslant \frac{1}{\left [ G:H \cap K \right ]}\underline{V_{\sigma}}\left ( H \cap K \right )\geqslant \frac{1}{\left [ G:H \right ]} V_{\sigma}\left ( H \right )}$$  This completes the proof of the proposition. \end{proof}

\vspace{0.2cm}
 
\begin{rems}${\textup{(i)}}$ By ${[12, \textup{Theorem 3.5}]}$, the mapping ${\omega}$ is well-defined and ${\omega\left ( G_{1} \right )=}$ ${\omega\left ( G_{2} \right )}$, if ${G_{1},G_{2}}$ are isomorphic groups in ${ \textup{V} \mathcal{F}}$. Also, the preceding remarks about the Kurosh-rank of a finite index subgroup of a free product show that ${\omega}$ is an Euler characteristic.

\vspace{0.2cm}

\noindent ${\textup{(ii)}}$ Let ${G}$ be a finitely generated group such that the net ${\left (\frac{\sigma\left ( H \right )}{\left [ G:H \right ]}\right )_{H}}$ is bounded and let ${\left (K_{n} \right )}$ be a decreasing and cofinal sequence of finite index subgroups of ${G}$. Then  we have the inequality ${V_{\sigma}\left ( G \right )\geqslant \varliminf_{n\rightarrow \infty}\frac{\sigma\left ( K_{n} \right )}{\left [ G:K_{n} \right ]}}$. Indeed, for every natural number ${n}$, there exists a finite index subgroup ${H_{n}}$ of ${G}$ such that ${V_{\sigma}\left ( G \right )+\frac{1}{n}\geqslant \sup_{K<H_{n}}\frac{\sigma\left ( K \right )}{\left [ G:K \right ]}}$. By the cofinality of ${\left (K_{n} \right )}$ we can construct inductively a subsequence with the property ${K_{s_{n+1}}\leqslant K_{s_{n}} \cap H_{n+1}}$ and hence ${V_{\sigma}\left ( G \right )+\frac{1}{n}\geqslant \varliminf_{n\rightarrow \infty} \frac{\sigma\left ( K_{s_{n}} \right )}{\left [ G:K_{s_{n}} \right ]}\geqslant \varliminf_{n\rightarrow \infty}\frac{\sigma\left ( K_{n} \right )}{\left [ G:K_{n} \right ]}}$.\\

\end{rems}

\vspace{0.2cm}

\noindent Let ${\chi}$ be an Euler characteristic defined on a class ${\mathcal{C}}$ with real values. Then, by Wall's argument we have the extension on ${\textup{V}\mathcal{C}}$: if ${G}$ is a group containing a finite index subgroup ${K}$ with ${K\in \mathcal{C}}$, then we define ${\chi\left ( G \right )=\frac{\chi\left (K \right )}{\left [ G:K \right ]}}$. The lemmas of the preceding section show that under certain conditions, the value of the extended mapping ${\chi}$ on amalgamated free products and HNN extensions over finite groups can be expressed in terms of the Euler characteristic of a free product with finitely many factors. In particular, if ${G}$ is residually finite or virtually torsion free, for any finite subgroup ${H}$ of ${G}$, there exists a normal and finite index subgroup ${K}$ such that ${K\cap H=1}$. Therefore, we have the following theorem:

\vspace{0.2cm} 

\begin{thm}
\label{thm:associativity}
Let ${\mathcal{C}}$ be a class of groups closed under free products and subgroups of finite index. We assume that ${\mathcal{C}}$ contains the infinite cyclic group ${\mathbb{Z}}$. Suppose that ${\chi:\mathcal{C}\rightarrow \mathbb{R}}$ is an Euler characteristic which satisfies the following property: ${\chi\left ( A \ast B\right )=\chi\left ( A\right )+\chi\left ( B\right )-1}$ for every ${A,B \in \mathcal{C}}$.

\vspace{0.2cm}

${\textup{(i)}}$ Let ${G\ast_{H}}$ be an HNN extension, where ${G}$ is residually finite or a virtually torsion free group in ${\mathcal{C}}$ and ${H}$ finite. Then ${G\ast_{H}}$ is virtually in ${\mathcal{C}}$ and: $${\chi\left ( G\ast_{H} \right )={\chi}\left ( G \right )- \frac{1}{|H|}}$$

\vspace{0.1cm}
${\textup{(ii)}}$ Let ${G_{1}\ast_{H} G_{2}}$ be an amalgamated free product, where ${G_{1}, G_{2}}$ are residually finite or virtually torsion free groups in ${{\mathcal{C}}}$ and ${H}$ finite. Then ${G_{1}\ast_{H} G_{2}}$ is virtually in ${\mathcal{C}}$ and: $${\chi\left ( G_{1} \ast_{H} G_{2} \right )={\chi}\left ( G_{1} \right )+{\chi}\left ( G_{2} \right )- \frac{1}{|H|}}$$
\end{thm}

\vspace{0.2cm} 

\noindent By Lemma 2.2, the class ${\textup{V} \mathcal{F}}$ is closed under free products and hence for the mapping ${\omega}$ we have the following result:
 
\vspace{0.2cm} 

\begin{cor} Let ${G}$ be the fundamental group of a finite graph of groups ${\left ( \mathcal{G}, Y \right )}$ with vertex groups in ${\textup{V} \mathcal{F}}$ and finite edge groups. We assume that each vertex group is either residually finite or virtually torsion free. Then the following formula holds: $${\omega\left ( G \right )=\sum _{u \in V Y}\omega\left ( G_{u} \right )-\sum _{e \in E Y}\frac{1}{\left | G_{e} \right |}}$$  In particular, if the  vertex groups of ${\left ( \mathcal{G}, Y \right )}$ are infinite, nilpotent and virtually torsion free, then isomorphic and finite index subgroups of ${G}$ have the same index.
\end{cor}

\vspace{0.1cm}

\begin{proof} Let ${G_{1},G_{2}}$ be groups which are virtually in ${\mathcal{F}}$. By definition, there exist finite index subgroups ${S_{1}=\ast_{i=1}^{n} A_{i}, S_{2}=\ast_{i=1}^{m} B_{i}}$ of ${G_{1},G_{2}}$ respectively such that the free factors ${A_{i},B_{j}}$ are in the class ${\mathcal{F}}$. Let ${d_{1}=\left [ G_{1}:S_{1} \right ]}$ and ${d_{2}=\left [ G_{2}:S_{2} \right ]}$. Then Lemma 2.2 ensures that the free product ${G_{1} \ast G_{2}}$ has a subgroup of index ${d_{1}d_{2}}$  which is a free product with ${\left ( n-1 \right )d_{2}+\left ( m-1 \right )d_{1}}$ ${+d_{1}d_{2}+1}$ free factors in ${\mathcal{F}}$. Then it is immediate that $${\omega\left ( G_{1} \ast G_{2} \right )=\frac{1-n}{d_{1}}+\frac{1-m}{d_{2}}-1=\omega\left ( G_{1} \right )+\omega\left ( G_{2} \right )-1}$$ Thus, the main theorem applies and the formula follows. \par For the second part, it is well-known that non-trivial nilpotent groups are freely indecomposable since they have non-trivial center, hence nilpotent groups are in ${\textup{V} \mathcal{F}}$ and ${\omega(G_{u})=0}$. Thus, by the previous part we have ${\omega\left ( G \right )=-\sum _{e \in EY}\frac{1}{\left | G_{e} \right |}<0}$. \end{proof}

\vspace{0.2cm}

\begin{lem}\emph{•} Let ${\mathcal{C}}$ be a class of groups closed under free products and subgroups of finite index, which contains the infinite cyclic group ${\mathbb{Z}}$. Suppose that ${\sigma:\mathcal{C}\rightarrow \mathbb{R}}$ is a mapping such that: for every ${G}$ in ${\mathcal{C}}$, the net ${\left ( \frac{\sigma\left ( H \right )}{\left [ G:H \right ]} \right )_{H}}$ is bounded and ${\sigma\left ( \mathbb{Z}\right )=0}$.
\vspace{0.2cm}

${\textup{(i)}}$ If ${\sigma\left ( A \ast B\right )\geqslant \sigma\left ( A\right )+\sigma\left ( B\right )+1}$ for every ${A,B \in \mathcal{C}}$ then $${\widetilde{V_{\sigma}}\left ( A \ast B \right )\geqslant \widetilde{V_{\sigma}}\left ( A\right )+\widetilde{V_{\sigma}}\left ( B\right )+1}$$

\vspace{0.2cm}

${\textup{(ii)}}$ If ${\sigma\left ( A \ast B\right )\leqslant \sigma\left ( A\right )+\sigma\left ( B\right )+1}$ for every ${A,B \in \mathcal{C}}$ then $${\widetilde{V_{\sigma}}\left ( A \ast B \right )\leqslant \widetilde{V_{\sigma}}\left ( A\right )+\widetilde{V_{\sigma}}\left ( B\right )+1}$$
\end{lem}

\vspace{0.1cm}

\begin{proof} There are three cases to consider.\\
\emph{Case 1:} We first assume that ${A,B}$ are infinite. Let ${\varepsilon >0}$. Then, there exist finite index subgroups ${H_{1}, H_{2}}$ of ${A,B}$ respectively, such that: ${\frac{\sigma\left ( H_{1} \right )}{d_{1}}>\underline{V_{\sigma}}\left ( A \right )}$ ${-\varepsilon}$, ${\frac{\sigma\left ( H_{2} \right )}{d_{2}}>\underline{V_{\sigma}}\left ( B \right )-\varepsilon}$, where ${d_{1}=\left [ A:H_{1} \right ],d_{2}=\left [ B:H_{2} \right ]}$. By Lemma 2.2, the free product ${A \ast B}$ has a subgroup ${K}$ of index ${d_{1}d_{2}}$ which is isomorphic to the free product ${\ast_{i=1}^{d_{2}}H_{1} \ast_{i=1}^{d_{1}}H_{2}}$ ${\ast F_{\left ( d_{1}-1 \right )\left ( d_{2}-1 \right )}}$ and thus the inequality ${\sigma\left ( A \ast B\right )\geqslant \sigma\left ( A\right )+\sigma\left ( B\right )+1}$ implies: $${1+\underline{V_{\sigma}}\left ( A \right )+ \underline{V_{\sigma}}\left ( A \right )-2\varepsilon\leqslant \frac{\sigma\left ( K \right )}{d_{1}d_{2}}\leqslant \underline{V_{\sigma}}\left ( A \ast B \right )}$$ Furthermore, there exists a finite index subgroup ${H}$ of ${A \ast B}$ such that ${\overline{V_{\sigma}}\left ( A \ast B \right )\geqslant}$ ${ \frac{\sigma\left ( H \right )}{d}-\varepsilon}$, ${d=\left [ A\ast B:H \right ]}$. By Kurosh subgroup Theorem, ${H}$ is a free product of the form $${\ast_{i \in I}\big( H \cap g_{i}Ag_{i}^{-1} \big) \ast \big(\ast_{i \in J}\big( H \cap w_{j}Aw_{j}^{-1} \big)\big)\ast F}$$ where ${g_{i}}$ ranges over a set of ${\left ( H,A \right )-}$double coset representatives and ${w_{j}}$ ranges over a set of ${\left ( H,B \right )-}$double coset representatives. We remark that the free factors of ${H}$ are non-trivial since ${A,B}$ are infinite, hence ${\textup{Kr}\left ( H \right )=|I|+|J|+\textup{rank}\left ( F \right )=1+d}$. Therefore, we have:
$${\overline{V_{\sigma}}\left ( A\ast B \right )+\varepsilon\geqslant}$$ $${\frac{1}{d}\sum_{i \in I}\sigma \big( H \cap g_{i}Ag_{i}^{-1}  \big)+\frac{1}{d}\sum_{j \in J}\sigma \big( H \cap w_{j}Bw_{j}^{-1} \big)+\frac{\textup{Kr}\left ( H \right )-1}{d}\geqslant }$$
$${\frac{\overline{V_{\sigma}}\left ( A \right )}{d}\sum_{i \in I}\big[ A: A \cap g_{i}^{-1}Hg_{i} \big]+\frac{\overline{V_{\sigma}}\left ( B \right )}{d}\sum_{j \in J} \big[ B: B \cap w_{j}^{-1}Hw_{j} \big]+1}$$ $${=\overline{V_{\sigma}}\left ( A \right )+\overline{V_{\sigma}}\left ( B \right )+1}$$

\vspace{0.1cm}

Thus, the inequality follows since ${\varepsilon}$ was arbitrary and ${V_{\sigma}\left ( G \right )=\inf_{H \in \Lambda_{G}}\left (\frac{1}{\left [ G:H \right ]} \underline{V_{\sigma}} \left ( H \right )\right )}$.\\

\noindent  \emph{Case 2:} ${A}$ is finite and ${B}$ is infinite. By Lemma 2.2, the free product ${A \ast B}$ has a subgroup of index ${|A|}$ which is isomorphic to ${\underbrace{B\ast...\ast B}_{|A| \ \textup{times}}}$, hence by the previous case we have that: ${V_{\sigma}\left ( A\ast B \right )=\frac{1}{\left | A \right |}V_{\sigma} \big( \underbrace{B\ast...\ast B}_{|A| \ \textup{times}} \big)\geqslant V_{\sigma}\left ( B \right )-\frac{1}{|A|}+1}$. If ${A}$ is infinite and ${B}$ is finite, we work similarly.\\

\noindent \emph{Case 3:} ${A}$ is finite and ${B}$ is finite. By Lemma 2.2, the free product ${A \ast B}$ has a free subgroup of index ${|A||B|}$ and rank ${r=\left ( |A|-1 \right )\left (|B|-1  \right )}$. Therefore, ${V_{\sigma}\left ( A\ast B \right )=\frac{1}{|A||B|}V_{\sigma}\left ( F_{r} \right )\geqslant 1-\frac{1}{|A|}-\frac{1}{|B|}}$.\\

\par The proof of the second part follows by using the same arguments.  \end{proof}

\vspace{0.2cm} 

Next we see some applications of the main theorem to some classical Euler characteristics.

\vspace{0.2cm}

Let ${G}$ be a finitely generated group. Lackenby in ${[4]}$ introduced the \emph{rank gradient} of ${G}$ to be $${\textup{RG}\left ( G \right )= \inf_{H \in \Lambda_{G}}\frac{r\left ( H \right )-1}{\left [ G:H \right ]}}$$\\ Since for every finite index subgroup ${H}$ of ${G}$ we have ${r\left ( H \right )-1\leqslant \left [ G:H \right ]\left ( r\left ( G \right )-1 \right )}$, it follows by Wall's observation ${\left [ \textup{12, Problem E10, p.85} \right ]}$ for sub-multiplicative invariants that the rank gradient is an Euler characteristic. The previous inequality also shows that ${\textup{RG}\left ( G \right )= \varlimsup_{H \in \Lambda_{G}}\frac{r\left ( H \right )-1}{\left [ G:H \right ]}}$.\\ 

Therefore, Gruschko theorem, Theorem 3.5 and Lemma 3.7 immediately imply the following corollary:

\vspace{0.2cm}

\begin{cor}\emph{•}Let ${G}$ be the fundamental group of a finite graph of groups ${\left ( \mathcal{G}, Y \right )}$ with finitely generated vertex groups and finite edge groups. We assume that the each vertex group is either residually finite or virtually torsion free. Then the following equality holds: $${\textup{RG}\left ( G \right )=\sum _{u \in V Y}  \textup{RG}\left ( G_{u} \right )+\sum _{e \in E Y}\frac{1}{\left | G_{e} \right |}}$$ \end{cor}

\vspace{0.2cm}

\noindent We remark that similar formulas for the rank gradient with respect to an infinite descending chain of finite index subgroups have been proved by Pappas ${[10]}$, in the case where the vertex groups are residually finite and the edge groups are amenable.

\vspace{0.2cm}

Another Euler characteristic defined on the class of finitely generated groups is the first ${L^2-}$Betti number. For the definition and other properties of this invariant, we refer the reader to ${[6,9]}$. Let ${G=A\ast_{C}B}$ be the amalgamated free product of ${A,B}$ over a common subgroup ${C}$ with ${b_{1}^{\left ( 2 \right )}\left ( C \right )=0}$. L$\textup{\"u}$ck in ${[7]}$ proved the formula $${b_{1}^{\left ( 2 \right )}\left ( G \right )-\frac{1}{\left | G \right |}=b_{1}^{\left ( 2 \right )}\left ( A \right )+b_{1}^{\left ( 2 \right )}\left ( B \right )-\frac{1}{\left | A \right |}-\frac{1}{\left | B \right |}+\frac{1}{\left | C \right |}}$$ By using the previous formula and applying the first part of Theorem 3.5 for the Euler characteristic ${G \mapsto b_{1}^{\left(2 \right )}\left ( G \right )-\frac{1}{\left | G \right |}}$, we conclude the corresponding formula for HNN extensions ${A\ast_{C}}$ where ${A}$ is residually finite or virtually torsion free and ${C}$ is finite. Therefore:

\vspace{0.2cm} 

\begin{cor} \emph{•}Let ${G}$ be the fundamental group of a finite graph of groups ${\left (\mathcal{G}, Y \right )}$ with finitely generated vertex groups and finite edge groups. We assume that each vertex group is either residually finite or virtually torsion free. Then the following equality holds: $${b_{1}^{\left ( 2 \right )}\left ( G \right )-\frac{1}{\left|G \right|}=\sum _{u \in V Y}\left ( b_{1}^{\left ( 2 \right )}\left ( G_{u} \right )-\frac{1}{\left | G_{u} \right |} \right )+\sum _{e \in E Y}\frac{1}{\left | G_{e} \right |}}$$ \end{cor}

\vspace{0.2cm}

Following the definitions presented in the beginning of this section, we define the \emph{Betti volume} of a finitely generated group ${G}$ to be $${V_{b}\left ( G \right )=\varlimsup_{H \in \Lambda_{G}}\frac{b_{1}\left ( H \right )-1}{\left [ G:H \right ]}}$$\\ where ${b_{1} \left(G \right)}$ denotes the rank of the free part of the abelianization ${G/\left [ G,G \right ]}$. By Proposition 3.3, ${V_{b}}$ is an Euler characteristic and we have ${V_{b}\left ( G \right )\leqslant \textup{RG}\left ( G \right )<\infty}$. Clearly, for finitely generated groups ${G_{1},G_{2}}$ we have ${b_{1}\left (G_{1}\ast G_{2} \right )=b_{1}\left ( G_{1} \right )+b_{1}\left ( G_{2} \right )}$, hence by using  Theorem 3.5 and Lemma 3.7, we obtain the following corollary:

\vspace{0.2cm}

\begin{cor}\emph{•}Let ${G}$ be the fundamental group of a finite graph of groups ${\left ( \mathcal{G}, Y \right )}$ with finitely generated vertex groups and finite edge groups. We assume that each vertex group is either residually finite or virtually torsion free. Then the following equality holds: $${V_{b}\left ( G \right )= \sum _{u \in V Y}  V_{b}\left ( G_{u} \right )+\sum _{e \in E Y}\frac{1}{\left | G_{e} \right |}}$$ \end{cor}

\vspace{0.2cm}

For a non-trivial finite group ${C}$ and a finitely generated group ${G}$, we denote by ${\textup{Hom}\left ( G,C \right )}$ the set of all homomorphisms from ${G}$ to ${C}$. Clearly, we have the inequality ${\left |\textup{Hom}\left ( G,C \right ) \right |\leqslant |C|^{r\left ( G \right )}}$ since each homomorphism from ${G}$ to ${C}$ is determined by the images of the elements of a generating set of ${G}$. For an infinite and finitely generated group ${G}$, we define ${V_{C}\left ( G \right )}$ to be the upper limit  $${V_{C}\left ( G \right )=\varlimsup_{H \in \Lambda_{G}}\frac{\frac{1}{\log|C|}\log\left | \textup{Hom}\left ( H,C \right ) \right |-1}{\left [ G:H \right ]}}$$ \par The previous remarks show that ${V_{C}\left ( G \right )}$ is dominated by ${r\left ( G \right )-1}$ and by Proposition 3.3 is an Euler characteristic. Note that if ${G}$ has a sequence of finite index subgroups ${\left(K_{n} \right)}$ with ${\left [ G:K_{n} \right ]\rightarrow \infty}$, then ${V_{C}\left ( G \right )\geqslant 0}$.
 \par
Moreover, the universal property of the free product implies $${\left | \textup{Hom}\left ( G_{1}\ast G_{2},C \right ) \right |=\left | \textup{Hom}\left ( G_{1},C \right ) \right |\cdot \left | \textup{Hom}\left ( G_{2},C \right ) \right |}$$ Thus, we conclude the following corollary:

\vspace{0.2cm}

\begin{cor}\emph{•}Let ${C}$ be a (non-trivial) finite group and ${G}$ the fundamental group of a finite graph of groups ${\left ( \mathcal{G}, Y \right )}$ with finitely generated vertex groups and finite edge groups. We assume that each vertex group is either residually finite or virtually torsion free. Then the following equality holds: $${\widetilde{V_{C}}\left ( G \right )= \sum _{u \in V Y}  \widetilde{V_{C}}\left ( G_{u} \right )+\sum _{e \in E Y}\frac{1}{\left | G_{e} \right |}}$$ where ${\widetilde{V_{C}}\left ( G \right )=\left\{\begin{matrix}
V_{C}\left ( G \right ),\ \textup{if} \ G \ \textup{is infinite} \ \ \ \ \\ 
 -\frac{1}{\left | G \right |},\ \textup{if} \ G \ \textup{is finite} \ \ \ \ \ \ 
\end{matrix}\right.}$ \end{cor}

\vspace{0.2cm}

\begin{prop}\emph{•}\\ \noindent ${\textup{(i)}}$ Let ${G}$ be a residually finite group with infinitely many ends. Then isomorphic and finite index subgroups of ${G}$ have the same index. Moreover, every monomorphism ${\varphi:G\rightarrow G}$ is either an automorphism or the image ${\varphi \left(G \right)}$ has infinite index in ${G}$.

\vspace{0.2cm}

\noindent${\textup{(ii)}}$ Let ${C}$ be a finite group. Then there exists a sequence of finite index subgroups ${\left(\Gamma_{n} \right)}$ of ${G}$ such that ${\left | \textup{Hom}\left ( \Gamma_{n},C \right ) \right |}$ grows at least as ${\left | C \right |^{\varepsilon\left [ G:\Gamma_{n} \right ]}}$ for some  ${\varepsilon>0}$.
\end{prop}

\vspace{0.1cm}

\begin{proof}${\textup{(i)}}$ By Stallings' theorem about ends of groups, ${G}$ splits non trivially over a finite subgroup, i.e. ${G}$ splits as an amalgamated free product ${G_{1}\ast_{F} G_{2}}$ or as an HNN extension ${H_{1}\ast_{F}}$. Since ${G}$ has infinitely many ends, in the first case we have ${\left [ G_{1}:F \right ]\geqslant 3}$ or ${\left [ G_{2}:F \right ]\geqslant 3}$ and in the second ${\left [ H_{1}:F \right ]\geqslant 2}$. By Corollary 3.8, we have that ${\textup{RG}\left ( G_{1} \ast_{F} G_{2} \right )\geqslant \frac{1}{6\left | F \right |}}$, ${\textup{RG}\left ( H_{1} \ast_{F}  \right )\geqslant \frac{1}{2|F|}}$ and therefore in each case the rank gradient of ${G}$ is at least ${\frac{1}{6|F|}}$. Clearly, the proof applies to virtually torsion free groups with infinitely many ends.

\vspace{0.2cm}

\noindent ${\textup{(ii)}}$ By Corollary 3.11 and the argument of the first part we have that ${V_{C}\left ( G \right )>0}$. Let ${M_{n}}$ be the intersection of finite index subgroups of ${G}$ of index at most ${n}$. By definition, for every ${n \in \mathbb{N}}$ there exists ${\Gamma_{n}<M_{n}}$ such that ${\frac{\frac{1}{\log\left | C \right |}\log\left | \textup{Hom}\left ( \Gamma_{n},C \right ) \right |-1}{\left [ G:H_{n} \right ]}+\frac{1}{n}\geqslant \sup_{K<\Gamma_{n}}\frac{\frac{1}{\log\left | C \right |}\log\left | \textup{Hom}\left ( K,C \right ) \right |-1}{\left [ G:K \right ]}}$ ${\geqslant V_{C}\left ( G \right )>0}$.   \end{proof}
\vspace{0.2cm}

We close this section with some examples.

\vspace{0.2cm}

\begin{exmp} The fundamental group of the closed and orientable surface ${S_{g}}$ of genus ${g\geqslant 1}$ is the finitely presented group ${\left \langle a_{1},b_{1},...,a_{g},b_{g}|\left [ a_{1},b_{1} \right ]...\left [ a_{g},b_{g} \right ] \right \rangle}$. If ${H}$ is a subgroup of ${\pi_{1}\left ( S_{g} \right )}$ of finite index, say ${d}$, then it is well known that ${H}$ is isomorphic to the fundamental group of ${ S_{h}}$ where ${h=1+d\left ( g-1 \right )}$. Hence by L$\textup{\"u}$ck's approximation theorem and the fact that ${r\left ( \pi_{1}\left ( S_{g} \right ) \right )=2g}$, we obtain  $${\textup{RG}\left ( \pi_{1}\left ( S_{g} \right ) \right )=b_{1}^{\left ( 2 \right )}\left ( \pi_{1}\left ( S_{g} \right ) \right )=2\left ( g-1 \right )}$$ (see also ${\left [ 9,\textup{Example 5.1.2} \right ]}$). \par More generally, let ${\left ( \mathcal{G},Y \right )}$ be a finite graph of groups with the following properties: each vertex group ${G_{u}}$ has a subgroup of index ${d_{u}}$ which is isomorphic to the fundamental group of ${S_{g_{u}}}$ (${g_{u}\geqslant 1}$) and each edge group is finite. If ${G}$ is the fundamental group of ${\left ( \mathcal{G},Y \right )}$, then: $${\textup{RG}\left(G  \right)=V_{b}\left(G \right)=b_{1}^{\left ( 2 \right )}\left ( G \right )=\sum_{u \in VY}\frac{2\left ( g_{u}-1 \right )}{d_{u}}+\sum_{e \in EY}\frac{1}{|G_{e}|}}$$ Moreover, since ${\pi_{1}\left ( S_{g} \right )}$ is freely indecomposable, we conclude that ${\omega\left ( G \right )=-\sum_{e \in EY}\frac{1}{|G_{e}|}}$.\end{exmp}

\vspace{0.2cm}

\begin{exmp} Let ${C}$ be a finite group and let ${\widehat{C}}$ be the set of all isomorphism classes of complex and irreducible representations of ${C}$. Mednykh in ${[8]}$ proved the following counting formula for the number of homomorphisms from ${\pi_{1}\left ( S_{g} \right )}$ to ${C}$ $${\left | \textup{Hom}\left ( \pi_{1}\left ( S_{g} \right ), C\right ) \right |=|C|^{2g-1}\sum _{V \in \widehat{C}}\left ( \frac{1}{\textup{dim}V} \right )^{2g-2}}$$ with ${g\in \mathbb{N} \cup \left \{ 0 \right \}}$. Let ${\left (H_{n}  \right )}$ be a sequence of finite index subgroups of ${\pi_{1}\left ( S_{g} \right )}$ such that ${\lim_{n\rightarrow \infty}\left [\pi_{1}\left ( S_{g} \right ):H_{n}  \right ]=\infty}$. In other words, there exist a sequence of integers ${(h_{n})}$ such that ${H_{n}}$ is isomorphic to ${\pi_{1}\left ( S_{h_{n}} \right )}$ and ${\lim_{n\rightarrow \infty}h_{n}=\infty}$. Then, Mednykh's formula gives \\$${\lim_{n\rightarrow \infty}\frac{\frac{1}{\log|C|}\log\left | \textup{Hom}\left ( H_{n},C \right ) \right |-1}{\left [ \pi_{1}\left ( S_{g} \right ):H_{n} \right ]}=}$$ $${2\left ( g-1 \right )+\lim_{n\rightarrow \infty}\frac{g-1}{\log|C|}\cdot \frac{1}{h_{n}-1}\log\left ( \sum_{V \in \widehat{C}} \left ( \frac{1}{\textup{dim}V} \right )^{2h_{n}-2}\right )=2\left ( g-1 \right )}$$\\ since there exists the trivial one dimensional representation of ${G}$. Therefore, by the second part of Remarks 3.4, we have ${V_{C}\left ( \pi_{1}\left ( S_{g} \right ) \right )=\textup{RG}\left ( \pi_{1}\left ( S_{g} \right ) \right )}$. In particular, if ${G}$ is the fundamental group of the graph of groups ${\left ( \mathcal{G},Y \right )}$ of the previous example, by Corollary 3.11 and the previous calculation we obtain ${V_{C}\left ( G \right )=\textup{RG}\left ( G \right )}$.  \end{exmp}

\vspace{0.2cm}

\begin{exmp} A polycyclic group is a solvable group which is Noetherian, i.e every subgroup is finitely generated. Let ${G}$ be the fundamental group of a finite graph of groups ${\left ( \mathcal{G}, Y \right )}$  with finite edge groups and vertex groups which are either polycyclic or finite. Polycyclic groups are amenable, finitely presented and residually finite, hence by ${\left [ 4, \textup{Theorem 1.2} \right ]}$ they have zero rank gradient. Thus, the previous corollaries imply $${\textup{RG}\left ( G \right )=V_{b}\left ( G \right )=b_{1}^{(2)}\left ( G \right )=-\sum _{u \in VY}\frac{1}{\left | G_{u} \right |}+\sum _{e \in EY}\frac{1}{\left | G_{e} \right |}}$$ where ${\frac{1}{\left | G_{u} \right |}=0}$ if ${G_{u}}$ is infinite.\end{exmp}

\section{Euler characteristics and accessibility}
In this section we see some applications of the lemmas of the previous section. We recall that an edge ${e}$ of a graph of groups ${\left ( \mathcal{G}, Y \right )}$ is called trivial if the initial ${\partial_{0}e}$ and the terminal ${\partial_{1}e}$ vertices of ${e}$ are distinct and we have either ${G_{e}\cong G_{\partial_{0}e}}$ or ${G_{e}\cong G_{\partial_{1}e}}$. If the edges are non-trivial, then ${\left ( \mathcal{G}, Y \right )}$ is called minimal.

\vspace{0.2cm}

For a group ${G}$ we set: ${\left \| G \right \|=\textup{sup} \big\{ \left | H \right |: H \space  \textup{\ is a finite subgroup of} \ G \big\}}$. Then we can prove the following: 

\begin{cor} 
${\textup{(i)}}$ Let ${G}$ be a finitely generated residually finite group with ${\left \| G \right \|< \infty}$. Then ${G}$ is accessible. In particular, if ${G}$ splits as a minimal graph of groups ${\left ( \mathcal{G}, Y \right )}$ with finite edge groups, then ${Y}$ has at most  ${\left \| G \right \| \left ( \frac{1}{2}+\textup{RG}(G) \right )}$ edges.
\vspace{0.2cm}

\noindent ${\textup{(ii)}}$ Let ${G}$ be a finitely generated virtually torsion free group. Then ${G}$ is accessible. In particular, if ${G}$ is not virtually free and ${H}$ is a finite index torsion free subgroup, then every splitting of ${G}$ as a minimal graph of groups with finite edge groups has at most ${\textup{RG}(H)}$ edges.  \end{cor}

\vspace{0.1cm}

\begin{proof} ${\textup{(i)}}$ We assume that ${G}$ is not accessible. By ${\left [5,\textup{Proposition 4} \right ]}$, for every natural number ${n}$, there exists a minimal graph of groups ${(\mathcal{G}_{n}, Y_{n})}$ with ${n}$ edges such that ${G \cong \pi_{1}(\mathcal{G}_{n},Y_{n} ,T_{n})}$ where ${T_{n}}$ denotes a fixed spanning tree of ${Y_{n}}$. Let ${u_{0}}$ be a vertex of ${Y_{n}}$. Then, by the minimality of ${\mathcal{G}_{n}}$, there exists a bijection ${i_{n}:ET_{n}\rightarrow VY_{n}\setminus\left \{ u_{0} \right \}}$ with the property: the edge group ${G_{e}}$ is embedded as a proper subgroup of ${G_{i_{n}(e)}}$. Furthermore, we observe: ${\frac{1}{|G_{e}|}+ \textup{RG}\left ( G_{i_{n}(e)} \right ) \geqslant \frac{1}{\left \| G \right \|}}$ since if ${G_{i_{n}(e)}}$ is infinite we have ${\textup{RG}(G_{u})\geqslant 0}$ as well as, if ${G_{i_{n}(e)}}$ is finite then ${\frac{1}{|G_{e}|}\geqslant \frac{2}{\left | G_{i_{n}(e)} \right |}\geqslant \frac{2}{\left \| G \right \|}}$. Therefore, by using Corollary 3.8, we obtain the inequality: $${\textup{RG}(G) = \sum_{u \in VY_{n}}\textup{RG}(G_{u})+ \sum_{e \in EY_{n}} \frac{1}{\left | G_{e} \right |}}$$ $${=\textup{RG}\left ( G_{u_{0}} \right )+ \sum_{e \in ET_{n}} \left ( \frac{1}{\left | G_{e} \right |}+ \textup{RG}\left ( G_{i_{n}(e)} \right )\right )+\sum _{e \in EY_{n}\setminus  ET_{n}}\frac{1}{|G_{e}|}}$$ $${\geqslant \textup{RG}(G_{u_{0}})+\frac{n}{\left \| G \right \|} \geqslant -\frac{1}{2}+\frac{n}{\left \| G \right \|}}$$\\ which gives the desired contradiction if we let ${n}$ tend to infinity. Finally, ${G}$ is accessible. 

\vspace{0.2cm}

\noindent${\textup{(ii)}}$ If ${H}$ is a torsion free subgroup of finite index, then ${\left \| G \right \|\leqslant \left [ G:H \right ]}$. Hence, by applying the previous argument we conclude that ${G}$ is accessible. Now if ${G}$ is not virtually free, for every minimal splitting ${\left ( \mathcal{G}, Y \right )}$ there exists a vertex ${u}$ such that the associated vertex group ${G_{u}}$ is infinite. Thus, ${\frac{\left | EY \right |}{\left [ G:H \right ]}+\textup{RG}\left ( G_{u} \right )\leqslant \textup{RG}\left ( G \right )}$.  
\end{proof}

\vspace{0.2cm}

\begin{rem} Linnell in ${[5]}$ used homological methods in order to obtain similar inequalities. More precisely, it is proved that for the fundamental group ${G}$ of a finite and minimal graph of groups ${\left ( \mathcal{G}, Y \right )}$, the rank of ${\mathbb{Q}\mathfrak{g}}$ (where ${\mathfrak{g}}$ denotes the augmentation ideal of ${G}$) as a ${\mathbb{Q}G-}$module is greater or equal than ${\frac{1}{2}+\frac{1}{2}\sum_{e \in EY}\frac{1}{\left | G_{e} \right |}}$.\end{rem}

\vspace{0.2cm}

By using ${[17,\textup{Theorem 6.12}]}$, we will see an application of the previous inequalities to the fixed point subgroups of automorphisms of fundamental groups. For the definition of the complexity of a group which acts on a simplicial tree and the maximum complexity of an accessible group, we refer the reader to ${\left [16,\textup{Definition 1.1} \right ]}$. In particular, the maximum complexity ${C_{\textup{max}}\left ( G \right )}$ of an accessible group ${G}$ which admits a minimal splitting ${\left ( \mathcal{G}, Y \right )}$ satisfies (by definition) the inequality ${C_{\max}\left ( G \right )\leqslant \left | EY \right |+1}$.

\vspace{0.2cm}

\begin{cor}\emph{•}Let ${G}$ be the fundamental group of a finite graph of groups ${\left ( \mathcal{G}, Y \right )}$ with Noetherian, virtually torsion free vertex groups and finite edge groups. Then for every automorphism ${f}$,  the fixed subgroup ${\textup{Fix}\left ( f \right )}$ is accessible. In particular, there exists a constant ${C=C\left ( G \right )}$ such that ${C_{\textup{max}}\left ( \textup{Fix}\left ( f \right ) \right )\leqslant C}$ for every automorphism ${f}$ of ${G}$.\end{cor}
\vspace{0.1cm}
\begin{proof} It is well known that Noetherian groups different from the infinite dihedral group are freely indecomposable, since they cannot contain a free subgroup of rank two. By Lemma 2.3 and Lemma 2.4, ${G}$ is virtually torsion free and there exists a finite index, characteristic and torsion free subgroup ${N}$. In particular, ${N}$ is a free product of the form ${\ast_{i=1}^{m}T_{i} \ast F\left ( X \right )}$ where each ${T_{i}}$ is a non-trivial, Noetherian and torsion-free subgroup of a conjugate of a vertex group ${G_{u}}$. Now let ${f \in \textup{Aut}\left (G \right ) }$. By ${[17,\textup{Theorem 6.12}]}$, the fixed point subgroup of the automorphism ${f|_{N}}$ inherits a splitting ${M_{1}\ast...\ast M_{d}\ast F\left ( S \right )}$ ${=N \cap \textup{Fix}\left ( f \right )}$ with ${d+\left | S \right |\leqslant m+\left | X \right |}$. Since ${M_{j}}$ is a subgroup of a conjugate some free factor ${T_{i}}$, the group ${N \cap \textup{Fix}\left (f \right )}$ is finitely generated. In particular, each factor ${M_{j}}$ is a group in ${\mathcal{F}}$ and ${-\omega \left ( N \cap \textup{Fix}\left ( f \right ) \right )=d+|S|-1 \leqslant}$ ${m+|X|-1=-\omega\left ( N \right )}$. It follows that ${\textup{Fix}\left ( f \right )}$ is  accessible since it is virtually torsion free and let ${\left ( \mathcal{G}_{f}, Y_{f} \right )}$ be a minimal splitting with maximum complexity. Then, by Corollary 3.6, ${-\omega\left ( \textup{Fix}\left ( f \right ) \right )=}$ ${-\sum_{u \in VY_{f}}\omega\left ( G_{u} \right )+\sum_{e \in EY_{f} }\frac{1}{\left | G_{e} \right |}}$ and similarly as in Corollary 4.1, we have the inequality: $${-\omega\left (\textup{Fix}\left ( f \right ) \right )\geqslant -\frac{1}{2}+\frac{\left | EY_{f} \right |}{\left \| \textup{Fix}\left ( f \right ) \right \|} \geqslant -\frac{1}{2}+ \frac{\left | EY_{f} \right |}{\left [ \textup{Fix}\left ( f \right ):N \cap \textup{Fix}\left ( f \right ) \right ]}}$$ which implies: $${C_{\textup{max}}\left ( \textup{Fix}\left ( f \right ) \right )\leqslant \left | EY_{f} \right |+1\leqslant \frac{\left [ \textup{Fix} \left ( f \right ): \textup{Fix}\left ( f \right ) \cap N \right ]}{2}-\omega \left (\textup{Fix}\left ( f \right ) \cap N   \right )+1 }$$ $${\ \ \ \ \ \ \leqslant \frac{\left [ G:N \right ]}{2}-\omega \left ( N \right )+1.}$$  \end{proof}

\vspace{0.45cm}

\noindent \textbf{Acknowledgements:} I would like to thank Mihalis Sykiotis for helpful suggestions and discussions. I also thank Denis Osin for useful comments on an earlier draft of this paper.

\vspace{0.2cm}

\bibliographystyle{unsrt}

\end{document}